\RequirePackage{fix-cm}

 \documentclass[showkeys,preprint,nofootinbib]{revtex4}          
\usepackage{graphicx}
\usepackage{amssymb}
\usepackage{amsbsy}
 \usepackage{amsmath}
\usepackage{subfigure}
\usepackage{color}
\usepackage{enumitem}
\usepackage{easybmat}
\usepackage{xypic}

\newcommand{\R}{\mathbb{R}}

\newcommand{\rd}{\mathrm{d}}

\usepackage{MnSymbol,wasysym}

\begin{document}

\title[]{From modelling of systems with constraints \\ to generalized geometry  and back to numerics
}

\author{Vladimir Salnikov} 
\email{vladimir.salnikov@univ-lr.fr}
\affiliation{
LaSIE  -- CNRS \& University of La Rochelle,
Av. Michel Cr\'epeau, 17042 La Rochelle Cedex 1, France}

\author{Aziz Hamdouni}
\email{aziz.hamdouni@univ-lr.fr}
\affiliation{
LaSIE -- University of La Rochelle, 
Av. Michel Cr\'epeau, 17042 La Rochelle Cedex 1, France}



\begin{abstract}
In this note we describe how some objects from generalized  geometry appear in the qualitative analysis and numerical simulation of mechanical systems.
In particular we discuss double vector bundles and Dirac structures.
It turns out that those objects can be naturally associated to systems with constraints -- we recall the mathematical construction 
in the context of so called implicit Lagrangian systems.  We  explain how they can be used to produce new numerical methods, that we  call \emph{Dirac integrators}.

On a test example of a simple pendulum in a gravity field we compare the Dirac integrators with classical explicit and implicit methods, we pay special attention to conservation of constrains.
Then, on a more advanced example of the Ziegler column we show that the choice of numerical methods can indeed affect the conclusions of qualitative analysis of the dynamics of mechanical systems.
We also tell why we think that Dirac integrators are appropriate for this kind of systems by explaining the relation with the notions of geometric degree of non-conservativity and kinematic structural stability.


\end{abstract}

\maketitle

\section{Introduction / motivation}
\label{sec:intro}

It is impossible to overestimate the importance of modelling for understanding natural mechanisms and improvement of industrial technologies. To see how much these techniques are appreciated nowadays, 
it is enough to cite the chemistry Nobel Prize 2013\footnote{This Nobel Prize was given to  Martin Karplus, Michael Levitt and Arieh Warshel ``For the development of multiscale models for complex chemical systems''.}: the work is related 
to the method of Molecular Dynamics -- an approach to modelling the condensed matter, 
popular since the eighties.\footnote{``Protein engineering ... represents the first
major step toward a more general capability for
molecular engineering which would allow us to
structure matter atom by atom.'' -- \cite{ulmer}}  
One can also mention some very concrete problems coming from mechanics: \emph{fluid--structure interaction}
in engineering (e.g. studying off-shore wind turbines, stability, resonances and fatigue of dock structures and 
bridges,\footnote{The famous Tacoma Narrows Bridge inaugurated in July 1940 and destroyed by wind in November 1940.} 
phenomena related to counters in petrol pipelines, the construction of thermal reactors, medical applications),
the development and implementation of \emph{constraints} in robotics; and this list is in no case exhaustive... 
One can even talk about \emph{open systems}, where considering the dissipation of energy
is a very subtle question. 

Even if nowadays it is easy to have access to serious computational resources, development of efficient 
methods remains a real challenge. Let us name a couple of examples in this context: \\
For molecular dynamics, one is interested in long trajectories. 
The computations carried out nowadays concern the dynamics on the timescale of picoseconds, 
while to study the phenomena like \emph{diffusion in porous media} 
one needs to be able to go to microseconds and even seconds.
To control the parameters of the system (like energy conservation), often some artificial ``correctors''
are added, like mechanical thermostats (Berendsen, Nos\'e--Hoover), 
even despite the fact that this can produce non-physical effects (see \cite{GSS, GS} and references therein). \\
In the fluid--structure interaction problem, the main difficulty is the size of data one needs to handle. 
On the one hand, one needs to take into account the geometry of the interface between the solid and the fluid, thus, 
introduce a very fine mesh. It is also important to have good discretization in time to capture the dynamics and 
especially deformations of the solid, that can influence  the spacial mesh as well.
On the other hand studied systems are usually very large, if one compares them with the scale of the interface and 
deformations. One thus needs to work with enormous mass of information, that results in algorithmic and technological 
issues. \\
For constraint systems appearing in robotics, the main concern is that already on the modelling phase some equality-type conditions appear -- the concept which is very vague from the numerical point of view. And we will see in this paper that 
``naive'' methods may result in completely non-physical behavior. 
In all the mentioned cases, as well as for many others, one needs reliable numerical schemes, preferably with reasonable
computational cost.

This paper is a part of a global project, the goal of which  is to develop mathematical and numerical tools, 
appropriate to these kind of problems. More precisely, we plan to study geometric structures that 
appear naturally for such problems, and numerical methods that ``respect'' these structures. And this approach has already given some fruitful results.

 For   \emph{conservative mechanical systems} the Hamiltonian/Lagrangian formalism gives a rather convenient
  framework for analysis of qualitative properties (like for example integrability or stability).
    From the numerical point of view, the integrators called \emph{symplectic} (\cite{verlet, yoshida}) exist
for some decades already -- they permit to control energy conservation.
The idea behind the method is to construct a discretization that will automatically preserve a 
symplectic structure\footnote{We will define the geometric objects needed to understand the results in the body of this paper. For all the others we mention, see Appendix \ref{bundles} for details.} -- this will guarantee the conservation of a Hamiltonian 
function that oscillates in a small neighborhood of the value of total energy of the system.
The phenomenon persists even for large time intervals, in contrast to other numerical schemes even of higher order.
The Lagrangian counterpart of this construction is related to integrators called \emph{variational} (\cite{west}).
For the example of molecular dynamics given above, popular methods are based on Verlet integration,  (e.g. ``velocity leapfrog''), which is symplectic. This is historically one of the first examples 
that has shown the use of studying geometric structures in the context of reliable numerical methods.
A lot of works have shown that numerical schemes that mimic physical properties of differential equations 
are more robust. The challenge is then to understand what mathematical structures one needs to consider 
while discretizing the equations to guarantee the conservation of physical quantities.
Numerical schemes constructed in this manner are often called 
\emph{geometric integrators}, since the mathematical structures behind come from differential or algebraic geometry.

In present days, there is a number of works related to methods preserving geometric structures, like 
symplectic or multi-symplectic forms\footnote{See Appendix \ref{bundles}.} (\cite{multisympl}), 
Lie symmetries,  
first integrals that result from the Noether's theorem.
A good example of application of geometric ideas in the context of numerics, is the Burgers' equation; 
coming from fluid mechanics, it is also used in gas dynamics and traffic. 
The equation is invariant under the Galilei group, but classical integration schemes 
destroy this invariance, this can create numerical artifacts, and thus a need 
to use very fine meshes to obtain reasonable solutions.  
Another option is to use Lie-symmetry based integrators (\cite{chhay}),  
that permits to remove this parasite solutions even with a coarser mesh. 
For this equation the  ``shock-waves'' phenomenon 
is also well-known: the profile of the wave becomes vertical and this poses problems for numerical integration.
The phenomenon can be well explained using \emph{contact geometry} and \emph{jet spaces} (\cite{rubtsov}),
 that also give a good phase space to use in computations. 
More generally, geometric integrators have been shown (\cite{chhay2, dina})
 to be more robust for PDEs, this is a very important 
result, saying that these schemes can be good candidates for simulation of large volumes mentioned 
 in the introduction.

The purpose of this paper is to apply this pattern to systems with constraints. The educated guess  (due to the 
works \cite{YoMa1, YoMa2}) is that the appropriate geometry is related to Dirac structures. We are sketching this construction 
in the next section about implicit Lagrangian systems. The section \ref{sec:discr} is devoted to discretization of the construction: we present it on a simple example and explain how to proceed  for more advanced ones.
In the application section (\ref{sec:numerics}) we do benchmarking tests of the implemented algorithms on a model problem of a simple pendulum. Then, we turn to more serious problems: the so called Ziegler column which is a model for a beam subject to external forces -- there we observe an essential difference in the results of simulation, depending on the choice of the numerical methods. We comment on this discrepancy in relation to other geometric concepts.

\newpage 
\section{Implicit Lagrangian systems}
As announced above, in this section we sketch the construction of Dirac structures naturally associated to systems with constraints.

\subsection{Dirac structures}
Dirac structures were introduced by T.~Courant \cite{courant} with the initial motivation coming from mechanics. 
The idea was to consider simultaneously velocities and momenta of a mechanical system, which as one knows are not independent. 
In the mathematical language, one would talk about the tangent bundle to a manifold $TM$ for velocities and the cotangent bundle 
$T^*M$ for momenta. Consider the direct sum of the two --  the so called Pontryagin bundle $E  = TM \oplus T^*M$, 
equipped with the following operations\footnote{This is an example of the so called Courant algebroid which turns out to be rather generic.} on pairs of its sections: \\
1) symmetric pairing: 
\begin{equation} \label{pairing}
 <v \oplus \eta, v' \oplus \eta'> = \iota_{v'}\eta + \iota_v\eta'  
 \end{equation}
2) Courant-Dorfman bracket: \\[-1em]
\begin{equation} \label{bracket}
   [v \oplus \eta, v' \oplus \eta'] =  
   [v,v']_{\text{Lie}} \oplus ({\cal L}_v \eta' - \iota_{v'} \rd \eta ).
\end{equation}
Here $\iota_v$ denotes the contraction of a vector field with a differential form, and ${\cal L}_v$ -- the Lie derivative along the vector field.

 An \emph{almost Dirac structure} ${\mathbb D}$ is a maximally isotropic (Lagrangian) subbundle  ${\mathbb D}$
  of $E$, i.e. a subbundle of $E$ on which the pairing (\ref{pairing}) vanishes identically, and which is of maximal rank equal to $\dim(M)$.
  If moreover the subbundle ${\mathbb D}$ is closed with respect to the bracket (\ref{bracket}), it is called a \emph{Dirac structure}\footnote{There is some confusion in literature, since people tend to omit the word ``almost'' in the first part of this definition.}.

The examples include a trivial one of
   ${\mathbb D} = TM$, and some more interesting ones of a graph of differential two form or of a bivector; thus Dirac geometry describes 
   uniformly symplectic and Poisson manifolds. 
The first condition related to (\ref{pairing}) is basically studying the linear algebra of the fibers over each point of $M$. The second one  
is more involved, and is sometimes called the integrability condition for ${\mathbb D}$. However, in this paper only the first condition will be relevant, 
we are thus using the \emph{almost} Dirac structures.

\subsection{Double bundles and systems with constraints}

In some sense, the definition of Dirac structures is a way of saying that velocities and momenta are related, but is was not explored in 
the original works of Courant. One of the reasons is that the geometric construction turned out to be very rich by itself. 
The other one is that the
good framework for mechanics is not that straightforward, typically one has to extend the configuration space somehow and then consider the subbundles.
For example, in the construction by H.~Yoshimura and J.~Marsden (\cite{YoMa1}) the bundles are over the complete phase space, i.e. 
$M = T^*Q$ -- the cotangent bundle to some configuration manifold $Q$. We will review the essential ingredients of the construction below, and exhibit the double bundles entering the play.

\textbf{1. Constraints.}

To describe a mechanical system with constraints, one gives some conditions that restrict its coordinates and velocities:
$$ 
  \varphi^a(q, \dot q) = 0, a = 1, \ldots, m
$$
Geometrically, this means that the dynamics takes place not on the whole tangent bundle but on some distribution of it
 $\Delta_Q \subset TQ$. It is convenient to view it as vector fields in the kernel of some set of differential forms $\psi^a, a = 1, \ldots, m$, i.e. at each point $q$ of the configuration space $Q$,
$\Delta_Q(q) = \{v \in T_qQ  \; | <\psi^a(q), v>= 0, \;  \forall a \}$.  

Consider the following diagram:
$$
\xymatrix{
   TT^*Q  \ar[dd] \ar@/^/ [rr]^{T\pi}    &&  \ar@/^/ [ll]^{pullback} TQ     \ar [dd] \\
  \\
 T^*Q   \ar [rr]^{\pi}&&Q  
 & 
} 
$$ 
Let $\pi$ be a canonical projection $\pi \; \colon \; T^*Q \to Q$ (map on the bottom), and $T\pi$ its tangent map: $T\pi\; \colon \; TT^*Q \to TQ$ (on the top).
Then $\Delta_Q$, a subbundle of $TQ$ (upper right corner) can be pulled back by $(T\pi)^*$.
Denote by $\Delta_{T^*Q} \subset TT^*Q$ its preimage, and $\Delta^0_{T^*Q} $ the annihilator of the preimage. 
If $(q, p)$ are coordinates on $T^*Q$,
locally 
$$
\Delta^0_{T^*Q}(q,p) = \{ \alpha_{(q,p)} \in T^*_{(q,p)}T^*Q \; | <\alpha_{(q,p)}, w_{(q,p)}> = 0, \; \forall w_{(q,p)} \in \Delta_{T^*Q}(q,p)  \}
$$

The canonical symplectic form $\Omega$ on $T^*Q$  defines\footnote{See Appendix \ref{bundles} for details} a mapping $\Omega^\flat \;\colon\; TT^*Q \to T^*T^*Q $. \\
An \emph{almost Dirac structure constructed from the constraint distribution} $\Delta_Q$ is a subbundle of $TT^*Q \oplus T^*T^*Q$ defined by: 
\begin{eqnarray} \label{DDeltaQ}
\mathbb{D}_{\Delta_Q} ((q,p)) = &\{& (w_{(q,p)}, \alpha_{(q,p)})  \in T_{(q,p)}T^*Q \times T^*_{(q,p)}T^*Q  \;|   \nonumber \\
 && w_{(q,p)} \in \Delta_{T^*Q}(q,p), 
\alpha_{(q,p)} - \Omega^\flat(q,p)w_{(q,p)} \in \Delta^0_{T^*Q}(q,p)  )    \}
\end{eqnarray} 

\textbf{2. Physics of the system.} 

Let $L \; \colon \; TQ \to \R$ be the Lagrangian. governing the system.  
 Its differential defines a mapping $\rd L \;\colon\; TQ \to T^*TQ$. 
For the coordinates $(q,v)$ on $TQ$, 
 locally it reads: \\ 
$\rd L \;\colon\; (q, v) \mapsto (q, v, \frac{\partial L}{\partial q}, \frac{\partial L}{\partial v})$. 

Recall again (\cite{tul} or see Appendix \ref{bundles}) that the double bundles are mapped to each other by symplectomorphisms:
$\Omega^\flat \;\colon\; TT^*Q \to T^*T^*Q$ and $\kappa \;\colon\; TT^*Q \to T^*TQ$. Then denote their composition by
$\gamma := \Omega^\flat \circ \kappa^{-1} \; \colon \; T^*TQ \to T^*T^*Q$. 
Define the  \emph{Dirac differential} ${\cal D} L := \gamma_Q \circ \rd L$. 
Locally it reads: ${\cal D} L \; \colon \; (q, v) \to (q, \frac{\partial L}{\partial v}, -\frac{\partial L}{\partial q}, v)$.

\textbf{3. Dynamics of the system.}

In the usual setting of ordinary differential equations, the dynamics of the system is recovered by integrating some vector field, in the constraint case the situation is a bit more intricate. 
 Consider a \emph{partial vector field} $X$, i.e. a mapping \\
$X \; \colon \;  \Delta_Q \oplus Leg(\Delta_Q) \subset TQ \oplus T^*Q \to TT^*Q$, \\
where $Leg(\Delta_Q)$ is the image of $\Delta_Q$ by the Legendre transform. 

An \emph{implicit Lagrangian system} is a triple $(L, \Delta_Q, X)$, s.t. $(X, {\cal D}L) \in \mathbb{D}_{\Delta_Q}$ (see eq. \ref{DDeltaQ}). \\
A \emph{solution} is a curve $(q(t), v(t), p(t)) \in TQ \oplus T^*Q$ integrating $X$.

The straightforward computation putting together the above definitions, shows that locally this translates into the following four conditions: 
  \begin{eqnarray} \label{impl_lagr}
    \dot q \in \Delta_{Q}, &\quad& p = \frac{\partial L}{\partial v} \\
    \dot q = v,   &\quad& \dot p - \frac{\partial L}{\partial q} \in \Delta^0(q) \nonumber
  \end{eqnarray}

Two remarks are in place here: \\
First, there is a clear abuse of notations -- 
a better geometric interpretation of $X$ is a vector field on $T^*Q$, where $(q(t), p(t))$ integrates it in the usual sense.  
The curve $v(t)$ appears because of the Dirac structure, namely the condition ${\cal D}L - \Omega^\flat X \in \Delta^0_{T^*Q}$,  
among others mentioned above, forces $\dot q = v$.
The other way around, it can be viewed as $X(q, v, p)$, where $p$ is given by the Legendre transform, and $v$ is in the constraint distribution. \\
Second, the last condition in (\ref{impl_lagr}) is the one that actually governs the dynamics in a non-trivial way. The left-hand-side of the inclusion is the  usual term from Euler--Lagrange equations. The right-hand-side is responsible for the constraints: the condition of being in the 
annihilator of something precisely means being a combination of its generators. 
And as we discussed above, the annihilator is generated by some one-forms $\psi^a$, when they are coming from the constraints 
$\psi^a = \rd \varphi^a$, one recovers 
$$
\dot p - \frac{\partial L}{\partial q} = \sum_a \lambda_a \rd\varphi^a,
$$
and recognizes immediately the Lagrange multipliers.

\section{Discrete version} \label{sec:discr}

At this point one may think that the whole formalism of the previous section is a complicated way to recover the classical story, and thus provides no ``added value''.  This is far from being true: the main message is that each operation described above admits a discrete version. 
This idea comes out very naturally and has been also explored by some followers of the Marsden's approach (see for example \cite{leok}). 
The goal of this section is to give the description of the final result using the minimal number of technicalities. 
The strategy that we have chosen to do this, is to consider a very explicit example of a simple pendulum, which however permits to 
explain the method. Moreover, we will see how to use some freedom of the approach to suggest potential improvements. 

\subsection{Example of an implicit Lagrangian system  \label{sec:pen}}
We start by applying the recipe from the previous section to the simple pendulum -- a mass point attached to a fixed point by a massless inextensible rod. 

\textbf{Description of $\Delta_Q$ and $\Delta_{T^*Q}$.}\\
In this case $Q=\R^2$, the constraint set is given by $\phi(x, y) := x^2 + y^2 - l^2 = 0$. The distribution is globally given by 
the vector field proportional (with a smooth non-zero coefficient depending on $x$ and $y$) to $\xi = y\frac{\partial}{\partial x} - x\frac{\partial}{\partial y}$, and it is obviously in the kernel of $\psi = \frac{1}{2}\rd\phi = x\rd x + y\rd y$. \\
This $\psi$ spans (with smooth coefficients) a subbundle of $T^*Q$ that we denote $\Delta^0_Q$ -- the annihilator of $\Delta_Q$.

As usual, the subbundle of $TT^*Q$ is denoted by  $\Delta_{T^*Q} = \{  v_{(q, p)} = (q, p, \dot q, \dot p) \quad | \quad \dot q \in \Delta_{Q} \}$, which here can be described globally. Its annihilator 
$\Delta^0_{T^*Q} \subset T^*T^*Q$ is as follows: \\ 
$\Delta^0_{T^*Q} = \{  \alpha_{(q,p)} =  (q, p, \alpha, w) \quad | \quad \alpha \in \Delta^0_Q, w = 0\}$.

\textbf{Lagrangian differential and Legendre transform.}\\
The Lagrangian  is 
\begin{equation} \label{lagr}
L = \frac{m}{2}(\dot x^2 + \dot y^2) - mgy. 
\end{equation}
The associated Lagrangian differential \\
${\cal D}L = ( q, \frac{\partial L}{\partial v}, -\frac{\partial L}{\partial q}, v ) = ((x, y), (m\dot x, m\dot y), (0, mg), (\dot x, \dot y))$.

\textbf{Almost Dirac structure.}\\
As described previously (eq. \ref{DDeltaQ}), we consider the almost Dirac structure $\mathbb{D}_{\Delta_Q}$ 
spanned by 
the couples $(v_{(q,p)}, \alpha_{(q,p)}) \in TT^*Q \oplus T^*T^*Q$ such that 
$\alpha_{(q,p)} - \Omega^\flat v_{(q,p)} \in \Delta^0_{T^*Q}$. 
Since locally $\Omega^\flat(v_{(q, p)}) = (q, p, -\dot p, \dot q))$, this condition can be rewritten explicitly as follows:
 at each point $(q, p)$,
$
   \mathbb{D}_{\Delta_Q}(q,p) = \{ ((q, p, \dot q, \dot p),(q, p, \alpha, w)) \quad | \quad \dot q \in 
   \Delta_Q(q, p), \; w= \dot q,  \; \alpha + \dot p \in \Delta^0_Q((q, p)) \}.
$

\textbf{All together.}\\
A vector field $X \in TT^*M$ defines the dynamics of the system in the sense of implicit Lagrangians when, together with 
${\cal D}L$ it belongs to  $\mathbb{D}_{\Delta_Q}$. 
Putting together the two previous paragraphs, we can rewrite the set of conditions (\ref{impl_lagr}) explicitly:
  \begin{eqnarray} \label{dirac_pen}
    \dot q \in \ker (\rd \varphi), &\quad& p = \frac{\partial L}{\partial v} \\
    \dot q = v,   &\quad& \dot p - \frac{\partial L}{\partial q} \approx \rd \varphi \nonumber
  \end{eqnarray}
The first one in \eqref{dirac_pen} precisely means that the velocity field should be compatible with $\varphi = 0$, which
is the differential view on the constraints. 
And the last one is again particularly interesting, since it says that $\dot p - \frac{\partial L}{\partial q} $ 
  should be proportional to the (only) generator of the constraint ideal $\psi = \frac{1}{2}\rd \phi$. In coordinates this latter reads:
\begin{eqnarray}
   \ddot x&=& \lambda x \nonumber \\
   \ddot y&=& - mg + \lambda y \nonumber
\end{eqnarray}
And we recognize again the Lagrange multiplier (here there is just one) in the right-hand-sides of the equations.

\subsection{Discretization  \label{subsec:discr}}
For the discretization, as mentioned above, we will give the final result and explain how it is obtained.
Let us recall the two imporatant continuous objects in the construction: \\
the  Lagrangian: $L = \frac{m}{2}(\dot x^2 + \dot y^2) - mgy, $ 
and the constraint: $ \phi(x, y) := x^2 + y^2 - l^2 = 0,$ giving rise to the differential condition 
$\iota_{\dot q} \psi \equiv \frac{1}{2} \iota_{\dot q} \rd \varphi = x \dot x + y \dot y = 0$
\\
Their straightforward discrete analogues are respectively:
$$
L_d (q_k, q_k^+)  = h\left(\frac{m}{2}\left(\frac{x_k^+ - x_k}{h}\right)^2 + \frac{m}{2}\left(\frac{y_k^+ - y_k}{h}\right)^2 - mgy_k\right)
$$ 
and
$
  \psi_d(q_k)  \equiv x_k \rd x + y_k \rd y = 0, 
  \qquad (\iota_{\dot q} \psi)_d(q_k, q_k^+) \equiv 
  x_k  \frac{x_k^+ - x_k}{h} + y_k \frac{y_k^+ - y_k}{h} = 0.
$\\
Let us explain the notations: the subscript ``$d$'' stands for ``discrete'', the index ``$k$'' represents the 
step number with the timestep $h$, the superscript ``$+$'' denotes the approximation of the corresponding variable for the next step.
The mathematical formulation of this result would be 
$L_d (q_k, q_k^+) = h L(q_k, {\cal R}_{q_k}^{-1}(q_k^+))$,  where 
${\cal R}$ is the retraction map from $TQ$ to $Q$, its inverse then is a way to reconstruct a tangent vector from two points on a manifold. The similar expression holds for the constraint.

The discretization of the equations \eqref{dirac_pen} or more generally \eqref{impl_lagr} is now rather natural: 
$\frac{\partial L}{\partial q}$ is replaced by $\frac{\partial L_d}{\partial q_k}$ and 
$\frac{\partial L}{\partial v}$ by $\frac{\partial L_d}{\partial q_k^+}$. Denoting the velocity approximation by 
$\tilde q$ one obtains the following system from \eqref{dirac_pen}: 
\begin{eqnarray} \label{discr_dirac}
&& x_k \tilde x_k +   y_k \tilde y_k = 0  \\
&& p^x_{k+1} = m\tilde x_k, \qquad p^y_{k+1} = m\tilde y_k, \nonumber \\
&& p^x_k - m \tilde x^k = \lambda x^k, \qquad p^y_k - m \tilde y^k -hmg = \lambda y^k.  \nonumber
\end{eqnarray}

It is important to note that at each time step there are five equations for five unknowns 
$x_{k+1}, p^x_{k+1}, y_{k+1}, p^y_{k+1}$ and $\lambda$, since 
$\tilde q_k$ are typically functions of known variables (indexed $k$) and $q_{k+1}$.
To convince oneself that it is not redundant, note that the first and the third line together give an expression for 
$\lambda$, which then can be plugged to the other equations. Then the $q_{k+1}$ are computed and
 $p_{k+1}$ are obtained from the second line - this is the core of the \emph{designed numerical method}.

The situation described above is absolutely generic, the simple example is given only for pedagogical 
purposes. In more details, one always  obtains the same number of equations and unknowns in the 
system similar to \eqref{discr_dirac}. They are typically linear for 
holonomic constraints. 
Non-linearity may come from two sources: the first line (constraints) and the retraction operator (i.e. the expression for $\tilde q$). 
And even being linear, the choice of $\tilde q$ is actually an important source of freedom for 
designing new numerical methods. 
In what follows, for illustration purposes we will consider two near at hand choices:  
$\tilde q_k = \frac{q_{k+1} - q_k}{h}$ (we label it Dirac-1) and $\tilde q^k = \frac{q_{k+1} - q_{k-1}}{2h}$ (Dirac-2).
 For the first choice one reproduces the results by \cite{leok}.

\section{Application -- numerical tests} \label{sec:numerics}

In this section we will present the results of simulations using various numerical methods. 
We start with a test example of a simple pendulum, observe some issued that can appear in the computations, and present 
what effects they can induce for more involved systems. 
For pedagogic reasons, we will ``push the limits'', that is do the simulation naively as if we had no a priori knowledge about the systems, and intentionally choose the parameters to clearly observe the effects.

\subsection{Pendulum}
\begin{figure}[ht]  
\centering
  \includegraphics[width=0.7\linewidth]{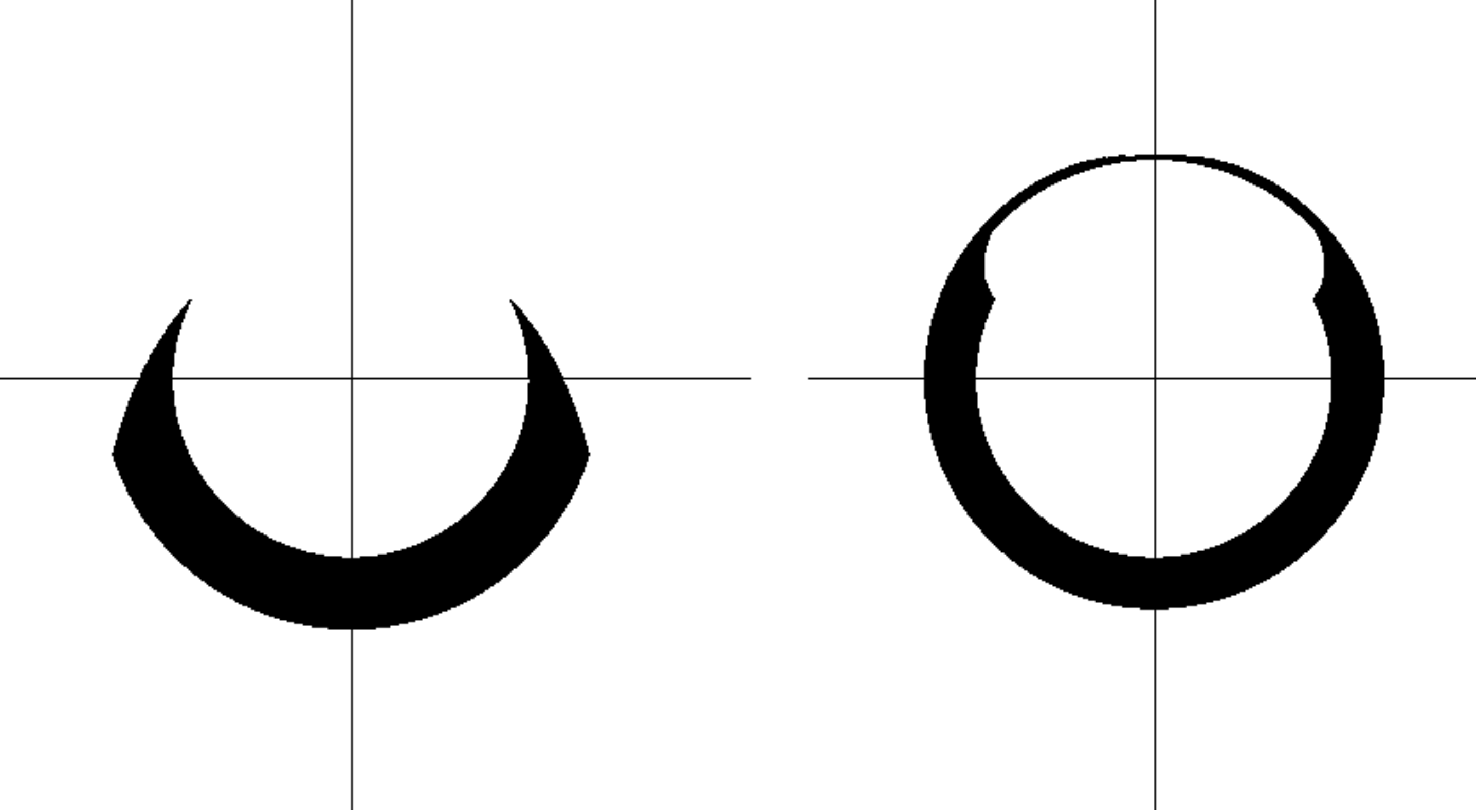}
\caption{\label{fig:pen-dirac1-euler} Region swept by the pendulum endpoint:  Dirac-1 (left) VS Euler (right).
} 
 \end{figure}

\begin{figure}[ht]  
\centering
  \includegraphics[width=0.7\linewidth]{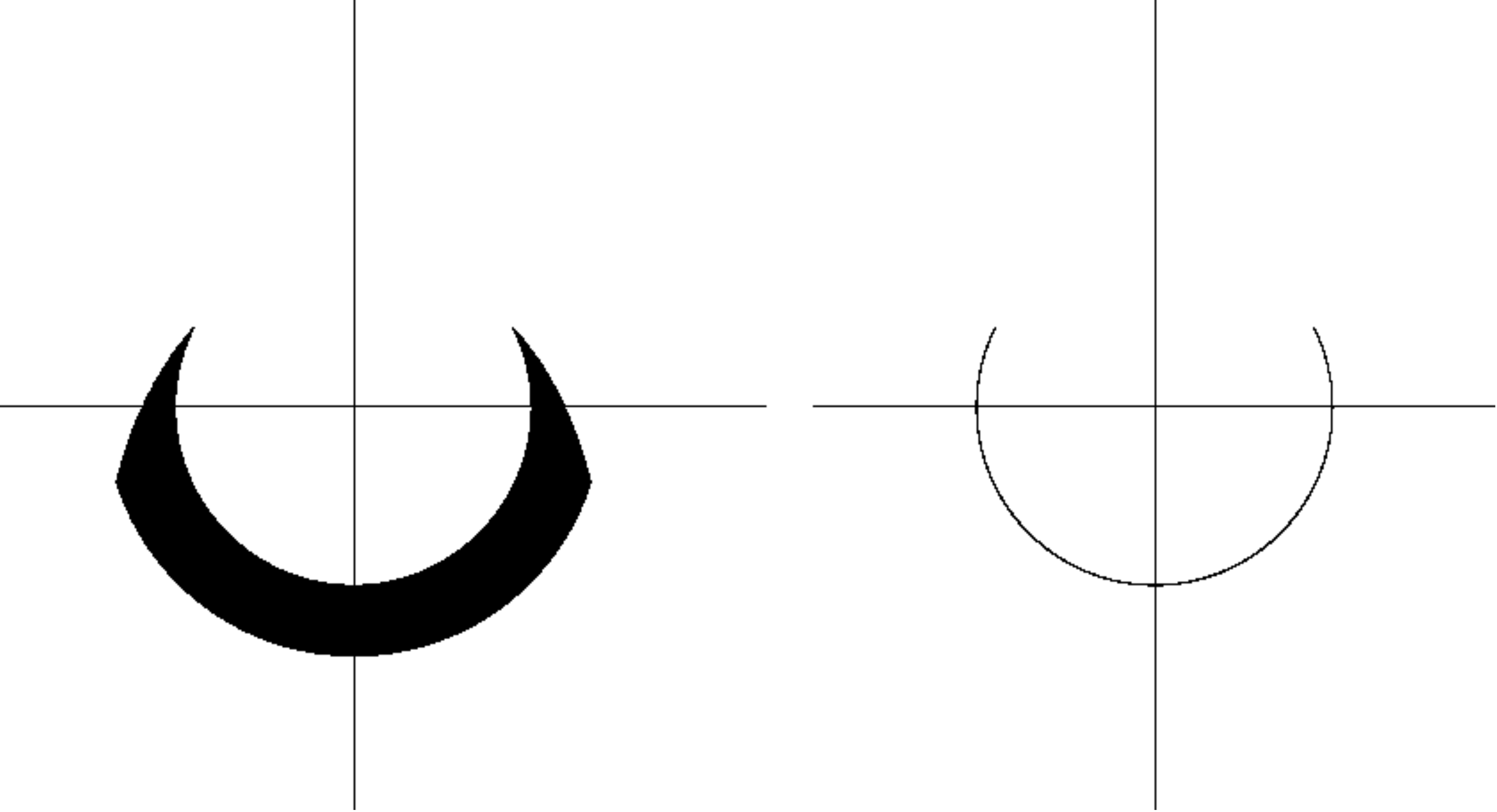}
\caption{\label{fig:pen-dirac1-dirac2} Region swept by the pendulum endpoint:  Dirac-1 (left) VS Dirac-2 (right).
} 
 \end{figure}

The system is the same as we studied in the section \ref{sec:pen} -- with a 2-dimensional configuration space, subject to 
one constraint. The trajectory thus should belong to the constraint level surface which is a curve.
We launch the simulation with the same initial conditions (zero velocity and some angle away from the equilibrium position) and the same timestep. The figures \ref{fig:pen-dirac1-euler} 
and \ref{fig:pen-dirac1-dirac2} show the regions swept by the trajectory after 100 periods of oscillation.
The  Dirac-1 method is compared respectively with the classical Euler method and the Dirac-2 method. 
One sees that the Euler method produces a totally non-physical behaviour: not only the constraints are violated, but also 
the dynamics is qualitatively different: the pendulum makes several full turns instead of oscillations. 
In this sense Dirac-1 is a bit better: at least the oscillatory nature of the dynamics is preserved. And Dirac-2 produces a much better picture -- the trajectory indeed sweeps a curve. 

Table \ref{tab} shows the relative constraint error for various numerical methods. It is important to note that although we are comparing the Dirac discretization procedure with classical methods, it is not completely honest. 
As described in subsection \ref{subsec:discr}, all the unknowns in the iteration procedure are treated in a similar way, 
while for traditional methods (Euler included), one needs to produce the explicit expression for Lagrange multipliers, which is often possible but not guaranteed. 

So, the preliminary conclusion here is that the Dirac structures permit to produce conceptually new numerical methods, that are very competitive. 

\begin{table}
\centering
\begin{tabular}{|c|c|}
\hline
Dirac order 1 & $0.952$ \\
\hline
Dirac order 2  &   $0.00204$  \\
\hline
Trapezium order 2 (implicit)& $0.0128$  \\
\hline
Adams--Bashforth order 3 &   $ 0.00014$\\
\hline
Runge--Kutta order 4    &    $ 9.4\cdot10^{-8}$ \\
\hline
\end{tabular}
\caption{Constraint error after 100 periods of oscillation.}
\label{tab}
\end{table}

\subsection{Ziegler column}
\begin{figure}[ht]  
\centering
  \includegraphics[width=0.48\linewidth]{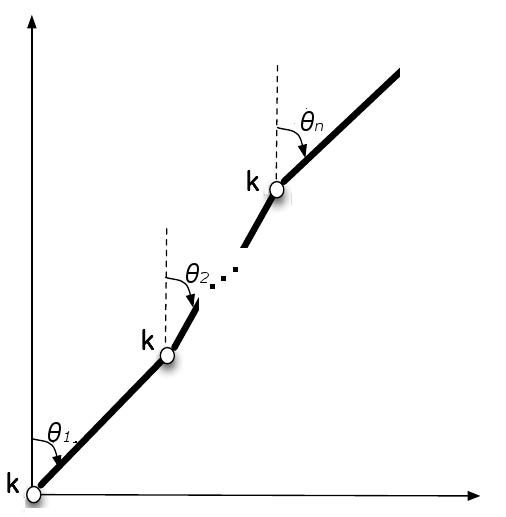}
\caption{\label{fig:zieg_3} Ziegler column.
} 
 \end{figure}

A Ziegler column (sometimes also called Ziegler--Bigoni system) is a chain of rods, attached to each other by endpoints (figure \ref{fig:zieg_3}). 
At each joint there is a harmonic force tending to align the segments, thus linear in the relative angle. External forces and torques 
may also be applied. We will consider it with a constant force $P$ applied to each joint, thus modelling a beam in a gravity field. 

Like in the previous section, we will start with an unconstraint system of $n$ mass points and fix the distances between the neighboring ones. 
Not to overload the presentation and also to simplify the visualization, we present the case $n=2$, the general one being absolutely similar. 

To apply the method from section \ref{sec:discr} we need the following data: the Lagrangian:
$$
L = \frac{m}{2}(\dot x_1^2 + \dot y_1^2) + \frac{m}{2}(\dot x_2^2 + \dot y_2^2) - P_c y_1 - P_c y_2 + \frac{k}{2} \theta_1^2 + \frac{k}{2} (\theta_2 - \theta_1)^2.
$$ 

where $\theta_i$ are functions of $x_1, x_2, y_2, y_2$; 
and the constraints: 
$$
\varphi_1 \equiv x_1 ^2 + y_1^2 - l_1^2 = 0, \quad \varphi_2 \equiv (x_2 - x_{2})^2 + (y_2 - y_{2})^2 - l_2^2 = 0. 
$$
the remaining part is merely repeating the discretization procedure from before. 

\begin{figure}[ht]  
\centering
  \includegraphics[width=0.7\linewidth]{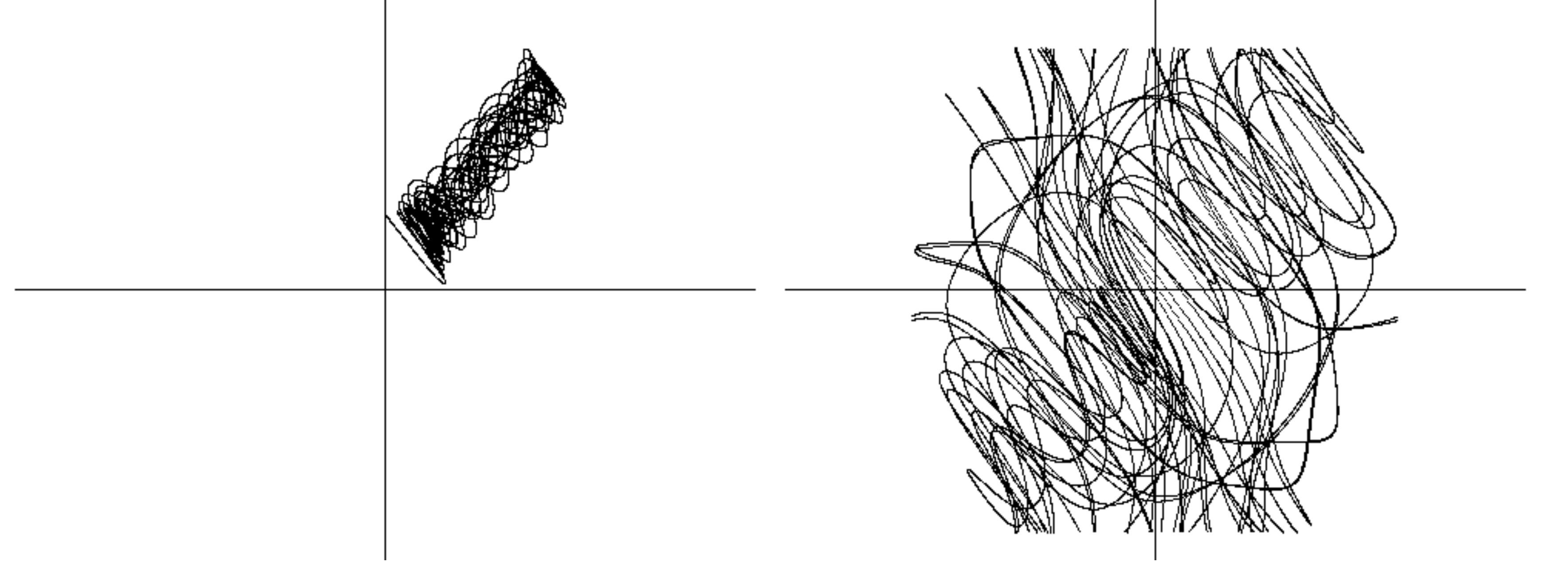}
\caption{\label{fig:zieg-dirac-trap} Evolution of angles of the system: Dirac (left) VS trapezium (right).
} 
 \end{figure}

\begin{figure}[ht]  
\centering
  \includegraphics[width=0.7\linewidth]{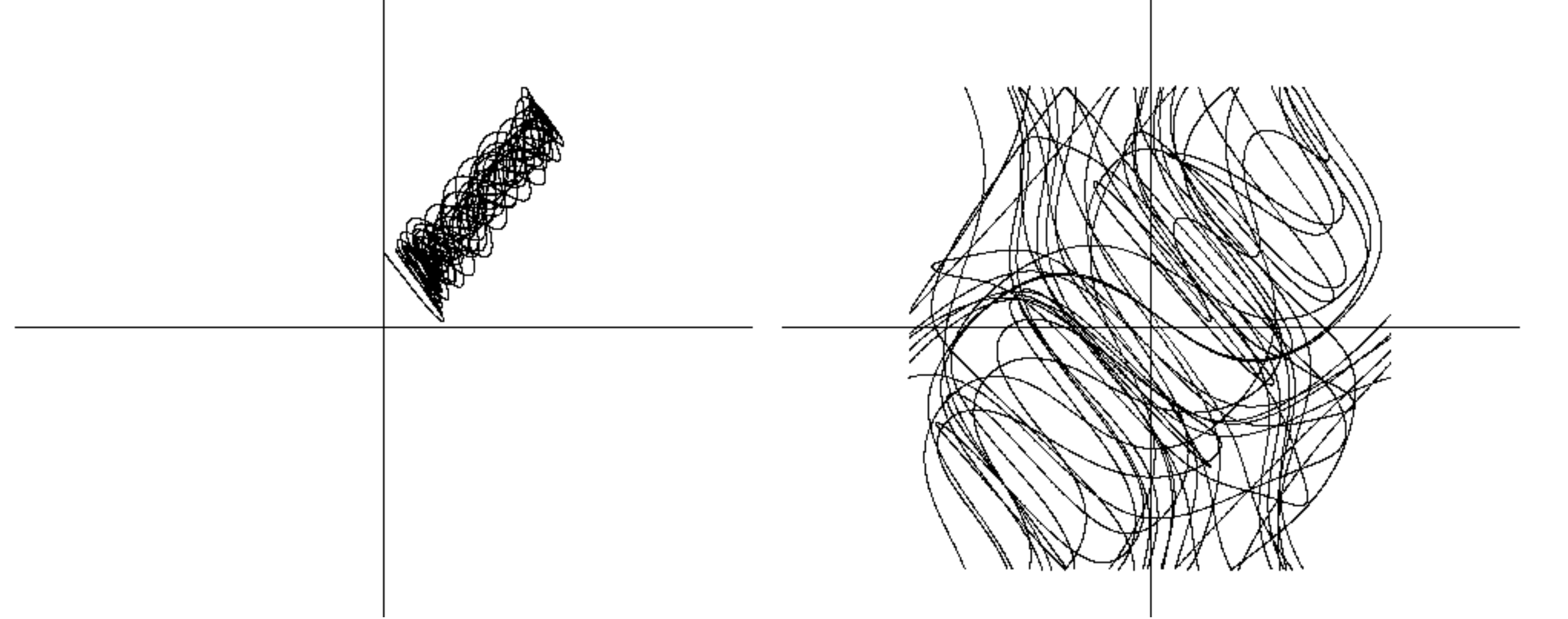}
\caption{\label{fig:zieg-dirac1-dirac2} Evolution of angles of the system: Dirac-2 (left) VS Dirac-1 (right).
} 
 \end{figure}

\begin{figure}[ht]  
\centering
  \includegraphics[width=0.7\linewidth]{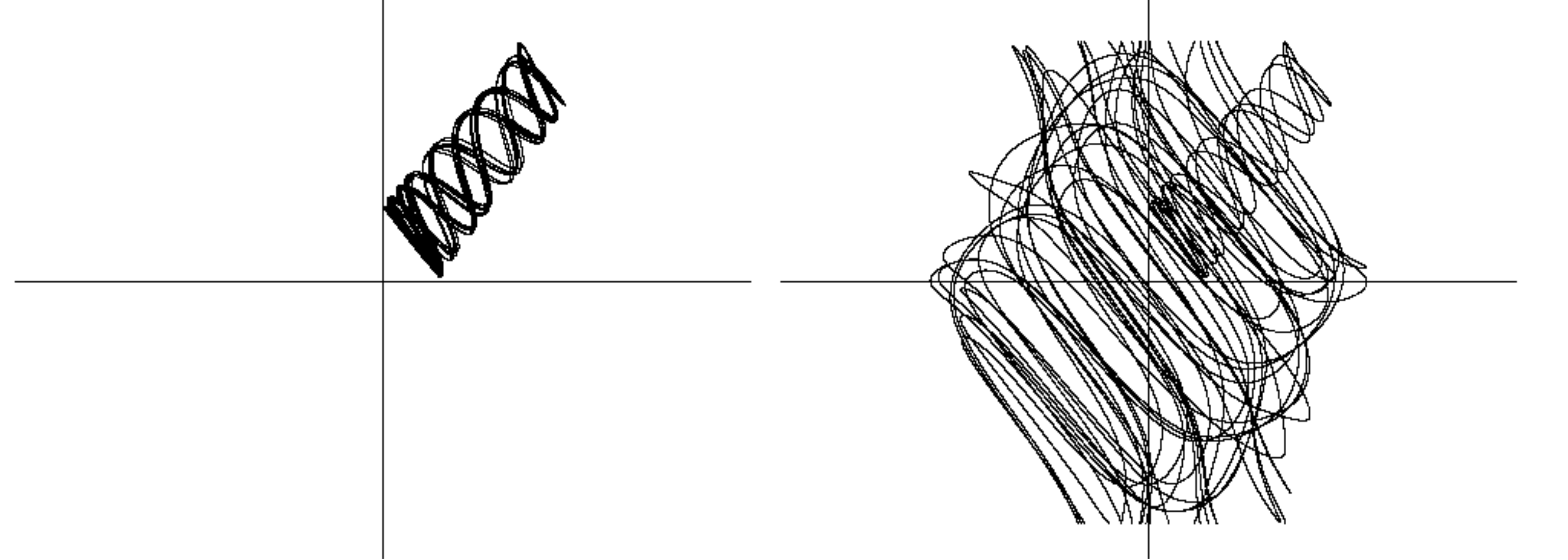}
\caption{\label{fig:zieg-gravity-stiffness}  Evolution of angles of the system: Dirac-2 (left) VS Dirac-1 (right).
} 
 \end{figure}

\begin{figure}[ht]  
\centering
  \includegraphics[width=0.7\linewidth]{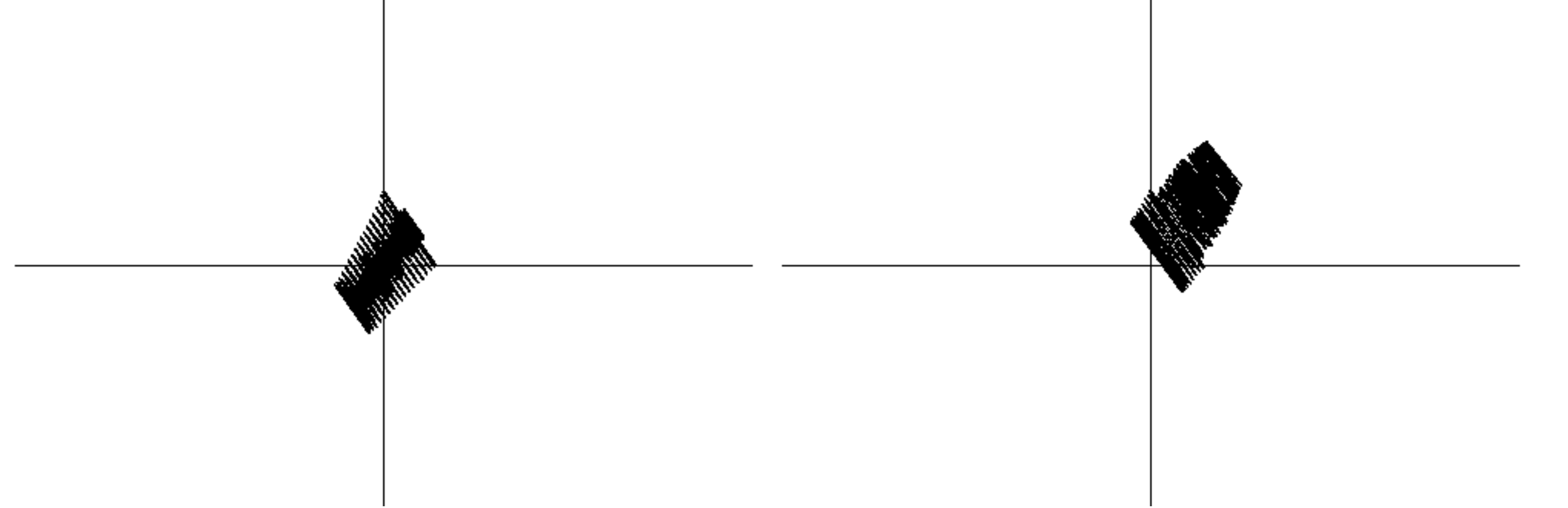}
\caption{\label{fig:zieg-stiffness-stiffness}  Evolution of angles of the system: Dirac-2 (left) VS Dirac-1 (right).
} 
 \end{figure}

Like for the pendulum, we will launch the simulations with the same initial data and parameters and compare the results. 
For visualization (figures \ref{fig:zieg-dirac-trap} -- \ref{fig:zieg-stiffness-stiffness}), we plot the angles computed using the coordinates of the endpoints.

In the dynamics of the system there are two competing forces: gravity and stiffness of the springs at the joints.
The first two figures (\ref{fig:zieg-dirac-trap} and \ref{fig:zieg-dirac1-dirac2}) illustrate basically the same phenomenon: 
one of the angles of the column reached an extreme value, and instead of bouncing back (correct behaviour on the left figure)
it made a full turn. This is clearly a numerical artifact that was produced precisely because of the accumulation of errors. 
The third figure (\ref{fig:zieg-gravity-stiffness}) is a similar effect, but the full turn happened not immediately. It would be fare to mention that these effects do not happen in real life, in a sense that even without the left part of the picture one would doubt the 
physics of the right one and dismiss the results of the computations. 

A more interesting effect happens on the figure 
\ref{fig:zieg-stiffness-stiffness} -- in both cases the system is oscillating, but around different equilibria. This is the most important illustration of the fact that the choice of a numerical method is crucial,  and that it is necessary to introduce the Dirac-2 method. The point is that there is no a priori reason to say that either of the solutions is non-physical, so in this case one is 
faced to the danger of producing a completely false solution without knowing it. 

Let us make another remark about the figure \ref{fig:zieg-stiffness-stiffness} here. In fact there was a series of works by 
J.~Lerbet and collaborators, studying analytically systems with constraints and in particular the Ziegler system (\cite{lerbet14} --\cite{lerbet18}). Various geometric concepts related to qualitative behaviour have been introduced there (Kinematic Structural Stability, Geometric Degree of Conservativity), and one can explain the observed effect in these terms. 
As we mentioned before, there are two competing forces: the stiffness, tending to stabilize the column in the upright vertical 
position, and the gravity, doing the opposite. And the balance depends on the ratio of these forces and the length of the constraints. It means that a slight change of the length due to the numerical error that we observed in section 
\ref{subsec:discr} can produce a bifurcation to a different equilibrium.

\section{Conclusion / discussions}

In this paper we have discussed a geometric approach to design and optimization of numerical methods for mechanical systems. 
We have seen that the naive application of classical methods may lead to obviously non-physical behaviour even for simple systems. 
The key point is as mentioned to encode the intrinsic physical properties via (differential) geometric structures and preserve the latter in the computation. The main advantage of this approach is that one controls more naturally the geometric structures, and more or less every discretization step has a geometric counterpart. 
For instance, for Dirac integrators there are some choices to make, like the approximating of the velocities related to the retraction map, and using this freedom one can improve the methods without guessing. 

Let us also mention that the usage of Dirac structures is not limited to systems with constraints described above.
We believe that this is a good language to unify the cited results by J.~Lerbet et al., for instance the geometric degree of conservativity is related to the rank of the projection of some well-chosen Dirac structure to the cotangent bundle. 
One important class 
where (almost) Dirac structures appear automatically, is the so called port-Hamiltonian systems (\cite{port-ham}). There the mathematics is well established, but numerics and benchmarking is rather intricate (\cite{dirac-port-ham}). 

Another activity in progress, inspired by \cite{dirac-port-ham} and the present paper is to
generalize the result: the goal is to define ``dynamical systems on Dirac structures'', i.e. some kind of higher analogue of Hamiltonian dynamics, where instead of symplectic structures the evolution is governed by the Dirac one.
To start with, an interesting question would be to relate this to the Poisson counterpart of the symplectic picture, and to the variational point of view of \cite{YoMa2}. 
And the ultimate goal would be to follow the same pattern ``physics $\to$ geometry $\to$ numerics'' for so called graded manifolds and $Q$-structures, which are known (see \cite{Sal, SalStr} and references therein) to give a unified description of symplectic, Poisson, and  Dirac structures.

\newpage
\appendix 
 \section{Symplectic structures and vector bundles }  \label{bundles}
For consistency, in this appendix we sketch some facts about double vector bundles and symplectic geometry in the context.
The results are mostly due to W.~Tulczyjew (\cite{tul}) and concern various iterations of tangent and cotangent bundles, that 
as we have seen appear naturally for the systems with constraints.
 
A \emph{symplectic manifold} is a manifold equipped with a \emph{symplectic form} -- a differential 2-form $\omega$ which is closed and non-degenerate. $\omega$ being closed means $\rd \omega = 0$, if moreover it is exact, i.e. the exists a 
1-form $\alpha$ such that $\omega = \rd \alpha$, this $\alpha$ is called the Liouville form.

For a symplectic 2-form,  $\omega$  the Hamiltonian vector field $X_f$ satisfies $\iota_{X_f} \omega = df$. 
This construction can be already useful for some PDEs like the stationary Lamb equation.
The generalization of this definition to differential $n+1$-forms
($n$-plectic) is $\iota_{X} \omega = d\alpha$, where
$\alpha$ is an ($n-1$)-form. This is an example of so-called \emph{higher} or \emph{generalized} geometry.

Important examples include the following construction.
Consider $\Omega = \rd p \wedge \rd q$ -- the canonical symplectic form on $T^*Q$, we omit the indeces for the coordinates, but all of them ($q, p, ...$) are multiplets of the dimension of the base manifold.
It induces the mapping 
$\Omega^\flat \colon TT^*Q \to T^*T^*Q$,
which in  coordinates reads: $(q, p, v_q, v_p) \mapsto (q, p, -v_p, v_q) = (q, p, \theta, \psi)$.  
On $T^*T^*Q$ there is the canonical symplectic form $\omega_1 = \rd \psi \wedge \rd p + \rd \theta \wedge \rd q$, 
with the Liouville form $\alpha_1 = \theta \rd q + \psi \rd p$. 
Hence, the induced symplectic structure on $TT^*Q$: $\tilde \omega_1 = \Omega^{\flat*}(\omega_1) = - \rd v_p \wedge \rd q + \rd v_q \wedge \rd p$ with the Liouville form $\tilde \alpha_1 = -v_p \rd q + v_q \rd p$.
Another natural mapping in the construction is the  Tulczyjew isomorphism 
 $\kappa \colon TT^*Q \to T^*TQ$.
In coordinates: $(q, p, v_q, v_p) \mapsto (q, v_q, v_p, p) = (q, v_q, \xi, \psi)$. 
On $T^*TQ$ there is the  canonical symplectic form as well: $\omega_2 = \rd \psi \wedge \rd v_q + \rd \xi \wedge \rd q$, 
with the Liouville form 
$\alpha_2 = \psi\rd v_q + \xi \rd q$.   
Hence, the induced symplectic structure on  $TT^*Q$ reads: $\tilde \omega_2 = \kappa^{*}(\omega_2) =  \rd v_p \wedge \rd q - \rd v_q \wedge \rd p$,
with the Liouville form $\tilde \alpha_2 = v_p \rd q + p \rd v_q$. 
Observe that $\tilde \alpha_1$ and $\tilde \alpha_2$ are essentially different, but $\tilde \omega_1 = - \tilde \omega_2$.

\newpage 
 

\end{document}